\documentclass[amscd,syntonly,amssymb]{amsart}

\usepackage{epsfig}
\theoremstyle{plain}
\newtheorem{Thm}{Theorem}
\newtheorem{Prop}[Thm]{Proposition}
\newtheorem{Cor}[Thm]{Corollary}

 \theoremstyle{definition}
\theoremstyle{remark}

\numberwithin{equation}{section}

\begin{document}
 \title{Hamiltonian diffeomorphisms of toric manifolds}

 \author{ ANDR\'{E}S   VI\~{N}A}
\address{Departamento de F\'{i}sica. Universidad de Oviedo.   Avda Calvo
 Sotelo.     33007 Oviedo. Spain. }
 \email{vina@uniovi.es}
\thanks{This work has been partially supported by Ministerio de Ciencia y
Tecnolog\'{\i}a, grant MAT2003-09243-C02-00}

 \maketitle
\begin{abstract}
We prove that $\pi_1(\text{Ham}(M))$ contains an infinite cyclic
subgroup, where $\text{Ham}(M)$  is the Hamiltonian group of the
one point blow up of  ${\Bbb C}P^3$. We give a sufficient
condition for the group $\pi_1(\text{Ham}(M))$ to contain an
infinite cyclic subgroup, when $M$ is a  general toric manifold.

\end{abstract}
\section {Introduction} \label{S:intro}

A  loop $\psi$  in the group $\text{Ham}(M,\omega)$ of Hamiltonian
symplectomorphisms \cite{Mc-S} of the symplectic manifold
$(M^{2n},\omega)$ determines a Hamiltonian  fibration
$E\stackrel{\pi}{\to}S^2$ with standard fibre $M$. On the total
space $E$ we can consider the first Chern class $c_1(VTE)$ of the
vertical  tangent bundle of $E$. Moreover on $E$ is also defined
the coupling class $c\in H^2(E,{\Bbb R})$ \cite{G-L-S}. This class
is determined by the following properties:

i) $i_p^*(c)$ is the cohomology class of the symplectic structure
on the fibre $\pi^{-1}(p)$, where $i_p$ is the inclusion of
$\pi^{-1}(p)$ in $E$ and $p$ is an arbitrary point of $S^2$.

ii) $c^{n+1}=0$.

These canonical cohomology classes determine the characteristic
number  \cite{L-M-P}
\begin{equation}\label{defIpsi}
I_{\psi}=\int_E c_1(VTE)c^n.
\end{equation}
 $I_{\psi}$  depends only on the homotopy class of $\psi$.
 Moreover
 $I$ is an ${\Bbb R}$-valued group homomorphism on
$\pi_1(\text{Ham}(M,\omega))$, so the non vanishing of $I$ implies
that the   group $\pi_1(\text{Ham}(M,\omega))$ is infinite. That
is, $I$ can be used to detect the infinitude  of the corresponding
homotopy group. Furthermore $I$ calibrates the Hofer's norm $\nu$
on $\pi_1(\text{Ham}(M,\omega))$ in the sense that $\nu(\psi)\geq
C|I_{\psi}|$, for all $\psi$, where $C$ is a positive constant
\cite{lP01}.

In \cite{V} we gave an explicit expression for the value of the
characteristic number $I_{\psi}$. This value can be calculated if
one has  a family of local symplectic trivializations of $TM$ at
one disposal, whose domains cover $M$ and are fixed by the
$\psi_t$'s. In fact we proved the  following Theorem
\begin{Thm}\label{Thefl}
 Let $\psi:S^1\to\text{Ham}(M,\omega)$ be a closed Hamiltonian isotopy.
If $\{B_1,\dots,B_m\}$ is a set of symplectic trivializations for
$TM$ which covers $M$ and such that $\psi_t(B_j)=B_j$, for all $t$
and all $j$, then
$$I_{\psi}=\sum_{i=1}^m J_i\int_{B_i\setminus\cup_{j<i}B_j}
\omega^n+\sum_{i<k}N_{ik},$$ where
$$N_{ik}=n\frac{i}{2\pi}\int_{S^1}dt\int_{A_{ik}}(f_t\circ\psi_t)
(d\, \log r_{ik})\wedge \omega^{n-1},$$
 $A_{ik}=(\partial
B_i\setminus\cup_{r<k}B_r)\cap B_k$, $J_{i}$ is the Maslov index
of $(\psi_t)_*$ in the trivialization $B_i$ and $r_{ik}$ the
corresponding transition function of $\text{det}\,(TM)$.
\end{Thm}

The homotopy type of $\text{Ham}(M,\omega)$ is completely known in
a few particular
 cases \cite{McD1} \cite{lP01} only.
 When $M$ is a surface,  $\text{Diff}_0(M)$ (the arc component of the identity map in
the diffeomorphism group of $M$) is homotopy equivalent to the
symplectomorphism group of $M$,  so the topology of the groups
$\text{Ham}(M)$ in dimension $2$ can be deduced from the
description of the diffeomorphism groups of  surfaces given in
\cite{E-E}
 (see \cite{lP01}). On the other hand,  positivity of
the intersections of
 $J$-holomorphic spheres in $4$-manifolds have been used in \cite{Gr} \cite{Ab} \cite{A-M} to prove
 results about the homotopy type of $\text{Ham}(M)$, when $M$ is a
 ruled surface. But these arguments
 which work in dimension  $2$ or dimension  $4$
 cannot be generalized to higher dimensions.

Let $\mathcal{O}$ be a coadjoint orbit of a Lie group $G$. If $G$
is semisimple and acts effectively on $\mathcal{O}$, McDuff and
Tolman have proved  that  the
 inclusion $G\to \text{Ham}(\mathcal{O})$ induces an injection from $\pi_1(G)$ to
 $\pi_1(\text{Ham}(\mathcal{O})) $ \cite{M-T}.
 This result answers a question posed in \cite{aW89}. In \cite{V1}
 we gave a lower bound for  $\sharp\pi_1(\text{Ham}(\mathcal{O}))$, when $\mathcal{O}$ is
a quantizable coadjoint orbit of a compact Lie group. In
particular we proved that $\sharp\pi_1(\text{Ham}({\Bbb
C}P^{n}))\geq n+1$.

In this note we use Theorem \ref{Thefl} to prove  that
 $\pi_1(\text{Ham}(M))$ contains an
 infinite cyclic subgroup, when $M$ is a particular toric  manifold.
 More precisely, when $M$ is the 6-manifold associated to the polytope
 obtained  truncating the tetrahedron of ${\Bbb R}^3$ with vertices
 $(0,0,0),\, (\tau,0,0),\, (0,\tau,0),\,  (0,0,\tau)$ by a horizontal plane \cite{Gui}; that
 is, when $M$ is the one point blow up of ${\Bbb C}P^3$. Moveover we will give
 a sufficient condition for $\pi_1(\text{Ham}(M))$ to contain an infinite cyclic
 subgroup, when $M$ is a general toric manifold.

The paper is organized as follows.  Section 2 is  concerned with
the determination of $I_{\psi}$ for a natural  circle action on
the one point blow up of ${\Bbb C}P^3$. In Section 3 we generalize
the arguments developed in Section 2 to toric manifolds. From this
generalization it follows
 the aforesaid sufficient
condition for the existence of an infinite subgroup in
$\pi_1(\text{Ham}(M))$, when $M$ is a toric manifold. Finally we
check that this sufficient condition does not hold for ${\Bbb
C}P^n$ with $n=1,2$. This  is consistent with the fact that
$\pi_1(\text{Ham}({\Bbb C}P^n))$ is finite for  $n=1,2$.

 I thank Dusa McDuff for her enlightening comments.


  \smallskip

\section{Hamiltonian group of the one point blow up of ${\Bbb C}P^3$}\label{blowup}
Given $\tau,\mu\in {\Bbb R}_{>0}$, with $\mu<\tau$, let $M$ be the
following manifold
  \begin{equation}\label{DefM}
  M=\{z\in{\Bbb
C}^5\,:\,|z_1|^2+|z_2|^2+|z_3|^2+|z_5|^2=\tau/\pi,\,\,|z_3|^2+|z_4|^2=\mu/\pi
\}/{\Bbb T}^2,
 \end{equation}
 where the action of ${\Bbb T}^2$ is
defined by
\begin{equation}\label{ActD}
 (a,b)(z_1,z_2,z_3,z_4,z_5)=(az_1,az_2,abz_3,bz_4,az_5),
 \end{equation}
  for
$a,b\in S^1$.

$M$ is a toric $6$-manifold; more precisely, it is the toric
manifold associated to the polytope obtained truncating the
tetrahedron of ${\Bbb R}^3$ with vertices
$$(0,0,0),\,(\tau,0,0),\,(0,\tau,0),\,(0,0,\tau)$$
  by a horizontal
plane through the point $(0,0,\lambda)$, with $\lambda:=\tau-\mu$
\cite{Gui}.

For $0\ne z_j\in{\Bbb C}$ we put $z_j=\rho_j e^{i\theta_j}$, with
$|z_j|=\rho_j$. On the set of points $[z]\in M$ with $z_i\ne 0$
for all $i$ one can consider the coordinates
\begin{equation}\label{Darb0}
\big(\frac{\rho_1^2}{2},\varphi_1, \frac{\rho_2^2}{2},\varphi_2,
\frac{\rho_3^3}{2},\varphi_3\big),
\end{equation}
 where the angle coordinates are defined by
  \begin{equation}\label{varph}
  \varphi_1=\theta_1-\theta_5,\,\varphi_2=\theta_2-\theta_5,\,
\varphi_3=\theta_3-\theta_4-\theta_5.
\end{equation}
Then the standard symplectic structure on ${\Bbb C}^5$ induces the
following form $\omega$ on this part of $M$
$$\omega=\sum_{j=1}^3d\big(\frac{\rho_j^2}{2}\big)\wedge d\varphi_j.$$

Let $0<\epsilon<<1$, we write
$$B_0=\{[z]\in M \,:\, |z_j|>\epsilon,\,\text{for all}\, j  \}.$$
For a given $j\in\{1,2,3,4,5\}$ we set
$$B_j=\{ [z]\in M \,:\,|z_j|<2\epsilon\;\,\text{and}\;\, |z_i|>\epsilon,\,\;\text{for all}\,\;i\ne j  \}$$
The family $B_0,\dots,B_5$ is not a covering of $M$, but if
$[z]\notin\cup B_k$, then there are $i,j$, with $i\ne j$ and
$|z_i|<\epsilon>|z_j|.$

On $B_0$ we will consider the well-defined Darboux coordinates
(\ref{Darb0}).
On $B_1$, $\rho_j\ne 0$ for $j\ne 1$; so the angle
coordinates $\varphi_2$ and $\varphi_3$ of (\ref{varph}) are
well-defined. We put $x_1+iy_1:=\rho_1 e^{i\varphi_1}$. In this
way we take as symplectic coordinates on $B_1$
$$\big(x_1,y_1, \frac{\rho_2^2}{2},\varphi_2,
\frac{\rho_3^3}{2},\varphi_3\big).$$
 We will also consider the following Darboux coordinates:
 On $B_2$
 $$ \big( \frac{\rho_1^2}{2},\varphi_1,x_2,y_2,
\frac{\rho_3^2}{2},\varphi_3\big),\;\text{with}\;
x_2+iy_2:=\rho_2e^{i\varphi_2}.$$
 On $B_3$
 $$ \big( \frac{\rho_1^2}{2},\varphi_1,
\frac{\rho_2^2}{2},\varphi_2,x_3,y_3\big),\;\text{where}\;
x_3+iy_3:=\rho_3e^{i\varphi_3}.$$
 On $B_4$
 $$ \big( \frac{\rho_1^2}{2},\varphi_1,
\frac{\rho_2^2}{2},\varphi_2, x_4,y_4\big),\;\text{with}\;
x_4+iy_4:=\rho_4e^{i\varphi_4}\;
\text{and}\;\varphi_4=\theta_4-\theta_3+\theta_5.$$
 On $B_5$
 $$ \big(x_5,y_5, \frac{\rho_2^2}{2},\chi_2,
\frac{\rho_3^2}{2},\chi_3\big),$$
 where
$$ x_5+iy_5:=\rho_5e^{i\varphi_5},\; \chi_2=\theta_2-\theta_1,\;
\chi_3=\theta_3-\theta_1-\theta_4,\;\chi_5=\theta_5-\theta_1.$$

If $[z_1,\dots,z_5]$ is a point of
 $$M\setminus\bigcup_{i=0}^5B_i,$$
  then there are $a\ne b\in \{1,\dots, 5\}$ such that $|z_a|,
  |z_b|<\epsilon$.
  We can cover the set $M\setminus\bigcup B_i$
  by Darboux charts denoted $B_6,\dots, B_q$ similar to the preceding
  $B_i$'s
satisfying the following condition:
  The image of each $B_a$, with $a=6,\dots, q$, is contained in a prism of ${\Bbb R}^6$ of the form
  $$\prod_{i=1}^6[c_i,\,d_i],$$
  where at least four intervals $[c_i,d_i]$ have length of order
  $\epsilon$.

  By the
infinitesimal ``size" of the $B_j$, for $j\geq 1$, it turns out
\begin{equation}\label{Oeps}
\int_{B_j}\omega^3=O(\epsilon),\;\text{for}\; j\geq 1.
\end{equation}

Let $\psi_t$ be the symplectomorphism  of $M$ defined by
\begin{equation}\label{DEfpsi}\psi_t[z]=[z_1e^{2\pi it},z_2,z_3,z_4,z_5].
 \end{equation}
  Then
$\{\psi_t\}_t$ is a loop in the group $\text{Ham}(M)$ of
Hamiltonian symplectomorphisms of $M$. By $f$ is denoted the
corresponding normalized Hamiltonian function. Hence
$f=\pi\rho_1^2-\kappa$ with $\kappa\in{\Bbb R}$ such that $\int_M
f\omega^3=0$.

In the coordinates (\ref{Darb0}) of $B_0$, $\psi_t$ is the map
$\varphi_1\mapsto\varphi_1+2\pi t$. So the Maslov index
$J_{B_0}=0$. It follows from (\ref{Oeps}) and Theorem \ref{Thefl}
\begin{equation}\label{IPN}
I_{\psi}=\sum_{i<k}N_{ik}+O(\epsilon),
\end{equation}
with
$$N_{ik}=\frac{3i}{2\pi}\int_{A_{ik}}fd\,\log r_{ik}\wedge
\omega^2.$$
 If $[z]\in A_{ik}\subset\partial B_i\cap B_k$ , with $1\leq i<k$,
 then at least the modules  $|z_a|$ and $|z_b|$ of two components of
 $[z]$ are of order $\epsilon$; so $N_{ik}$ is of order $\epsilon$
 when $1\leq i<k$.
Analogously $N_{0k}$ is of order $\epsilon$, for $k=6,\dots, q$.
Hence (\ref{IPN}) reduces to
\begin{equation}\label{IPNr}
I_{\psi}=\sum_{k=1}^5 N_{0k}+O(\epsilon).
\end{equation}
 If we put
\begin{equation}\label{N0kpri}
 N_{0k}'=\frac{3i}{2\pi}\int_{A_{0k}'}fd\,\log r_{ik}\wedge
\omega^2,
\end{equation}
 with
 $$A_{0k}'=\{[z]\in M\,:\, |z_k|=\epsilon,\, |z_r|>\epsilon \;\text{for all}\; r\ne k
 \}$$
then
\begin{equation}\label{auxi}
N_{0k}=N_{0k}'+O(\epsilon).
\end{equation}

Next we determine the value of $N'_{01}$. To know the transition
function $r_{01}$ one needs the Jacobian matrix $R$ of the
transformation
$$\big(x_1,y_1,\frac{\rho_2^2}{2},\varphi_2,\frac{\rho_3^2}{2},\varphi_3\big)\to
\big(\frac{\rho_1^2}{2},\varphi_1,\frac{\rho_2^2}{2},\varphi_2,,\frac{\rho_3^2}{2},\varphi_3\big)$$
in the points of $A'_{01}$; where $\rho_1^2=x_1^2+y_1^2$,
$\varphi_1=\tan^{-1}(y_1/x_1)$. The function $r_{01}=\rho(R)$,
where $\rho:Sp(6,{\Bbb R})\to U(1)$ is the map which restricts to
the determinant on $U(3)$  \cite{S-Z}.
 The non trivial block of $R$ is
the diagonal one
$$\begin{pmatrix} x_1& y_1 \\
r& s \\
\end{pmatrix},$$
 with $r={-y_1}{(x_1^2+y_1^2)}^{-1}$ and
$s={x_1}{(x_1^2+y_1^2)}^{-1}$. The non real eigenvalues of $R$ are
$$\lambda_{\pm}=\frac{x_1+s}{2}{\pm}\frac{i\sqrt{4-(s+x_1)^2}}{2}.$$
On $A'_{01}$  these non real eigenvalues occur when $(s+x_1)^2<2$,
that is, if $|\cos
\varphi_1|<2\epsilon(\epsilon^2+1)^{-1}=:\delta$. If $y_1>0$ then
$\lambda_{-}$ of the first kind (see \cite{S-Z}) and $\lambda_+$
is of the first kind, if $y_1<0$.

Hence, on $A'_{01}$,
\begin{equation}\notag \rho(R)=\begin{cases}
\lambda_+|\lambda_+|^{-1}=x+iy,
&\text{if $\,|\cos\varphi_1|<\delta\,$ and $\,y_1<0$;}\\
\lambda_-|\lambda_-|^{-1}=x-iy,
&\text{if $\,|\cos\varphi_1|<\delta\,$ and $\,y_1>0$;}\\
\pm 1,&\text{otherwise.}
\end{cases}
\end{equation}
where $x=\delta^{-1}\cos\varphi_1$, and $y=\sqrt{1-x^2}$.

If we put $\rho(R)=e^{i\gamma}$, then $\cos\gamma=\delta^{-1}\cos
\varphi_1$ (when $|\cos\varphi_1|<\delta$), and
\begin{equation}\notag
\sin\gamma=
\begin{cases}-\sqrt{1-\cos^2\gamma},&\text{if
$\,\sin\varphi_1>0$;}\\
\sqrt{1-\cos^2\gamma},&\text{if $\,\sin\varphi_1<0$.}
\end{cases}
\end{equation}
So when $\varphi_1$ runs anticlockwise from $0$ to $2\pi$,
$\gamma$ goes round clockwise the circumference; that is,
$\gamma=h(\varphi_1)$, where $h$ is a function such that
\begin{equation}\label{h1}
h(0)=2\pi,\;\; \text{and}\;\; h(2\pi)=0.
 \end{equation}
 As $r_{01}=\rho(R)$,
then $d\,log \,r_{01}=id h$.

On $A'_{01}$ the form   $\omega$ reduces to $(1/2)d\rho_2^2\wedge
d\varphi_2+d\rho_3^2\wedge d\varphi_3$. From (\ref{N0kpri}) one
deduces
\begin{equation}\label{N01a}
N'_{01}=\frac{3i}{4\pi}\int_{A'_{01}}i f \frac{\partial
h}{\partial \varphi_1}\,d\varphi_1\wedge d\rho_2^2\wedge
d\varphi_2 \wedge d\rho_3^2\wedge d\varphi_3 .
\end{equation}

 The submanifold $A_{01}'$ is oriented as a subset of $\partial B_0$ and the orientation of $B_0$
 is the one defined by $\omega^3$, that is, by
 $$d\rho_1^2\wedge d\varphi_1 \wedge d\rho_2^2\wedge
d\varphi_2 \wedge d\rho_3^2\wedge d\varphi_3.$$
  Since $\rho_1>\epsilon\,$  for the points of $B_0$, then $A_{01}'$
is oriented by $-d\varphi_1\wedge d\varphi_2^2\wedge
d\varphi_2\wedge d\rho_3^2\wedge d\varphi_3.$ On the other hand,
 the Hamiltonian function $f=-\kappa+O(\epsilon)$ on $A'_{01}$.
Then it follows from (\ref{N01a}) together with (\ref{h1})
$$N'_{01}=6\pi^2\kappa\int_{0}^{\mu/\pi}d\rho_3^2\int_{0}^{\tau/\pi-\rho_3^2}d\rho_2^2+O(\epsilon).$$
that is,
\begin{equation}\label{N01}
N'_{01}=3\kappa(\tau^2-\lambda^2)+O(\epsilon).
\end{equation}

The contributions $N'_{02}, N'_{03}, N'_{04}, N'_{05}$ to
$I_{\psi}$ can be calculated in a similar way. One obtains the
following results up to addends of order $\epsilon$
\begin{equation}\label{contr}
N'_{02}=N_{05}=-(\tau^3-\lambda^3)+3\kappa(\tau^2-\lambda^2)
,\;\;N'_{03}=\tau^2(3\kappa-\tau),\;\;N'_{04}=\lambda^2(3\kappa-\lambda).
\end{equation}
As $I_{\psi}$ is independent of $\epsilon$, it follows from
(\ref{IPNr}), (\ref{auxi}),  (\ref{N01}) and (\ref{contr})
\begin{equation}\label{Ipseo}
 I_{\psi}=6\kappa(2\tau^2-\lambda^2)+\lambda^3-3\tau^3.
\end{equation}

 On the other hand, straightforward  calculations give
$$\int_M \omega^3=(\tau^3-\lambda^3),\;\text{and}\;
\int_M\pi\rho_1^2\,\omega^3=\frac{1}{4}(\tau^4-\lambda^4).$$
 So
 \begin{equation}\label{kap}
\kappa=\frac{1}{4}\Big(\frac{\tau^4-\lambda^4}{\tau^3-\lambda^3}\Big).
\end{equation}
It follows from  (\ref{Ipseo}) and (\ref{kap})
\begin{equation}\label{IPSf}
I_{\psi}=\frac{\lambda^2(-3\tau^4+8\tau^3\lambda-6\tau^2\lambda^2
+\lambda^4)}{2(\tau^3-\lambda^3)}.
\end{equation}
Hence $I_{\psi}$ is a rational function of $\tau$ and $\lambda$.
It is easy to check that its numerator  does not vanish for
$0<\lambda<\tau.$ So we have
\begin{Prop}\label{Prononv}
If $\psi$ is the closed Hamiltonian isotopy defined in
(\ref{DEfpsi}), then  the   characteristic number $I_{\psi}\ne 0$.
 \end{Prop}

Next we consider the loop $\tilde \psi$ defined by
 \begin{equation}\label{Defetildepsi}
 \tilde\psi_t[z]=[z_1,z_2,z_3e^{2\pi it},z_4,z_5].
 \end{equation}
 The corresponding normalized Hamiltonian function is $\tilde
 f=\pi\rho_3^2-\tilde\kappa,$ where
 \begin{equation}\label{Kaptilde}
 \tilde\kappa=\frac{1}{4}\Big(\frac{\tau^4-4\tau\lambda^3+3\lambda^4}{\tau^3-\lambda^3}
 \Big).
 \end{equation}
 As in the preceding case
\begin{equation}\label{Ipsitil}
I_{\tilde\psi}=\sum_{j=1}^5 \tilde N'_{0j} +O(\epsilon),
\end{equation}
where
$$\tilde N'_{0j}=\frac{3i}{2\pi}\int_{A'_{0j}}\tilde f d\,\log
r_{0j}\wedge\omega^2.$$

The expression for $\tilde N'_{01}$ can be obtained from
  (\ref{N01a}) substituting  $f$ for $\tilde f$; so
  \begin{equation}\label{N01t}
  \tilde N'_{01}=-3(\tau-\tilde\kappa)(\tau^2-\lambda^2)+2(\tau^3-\lambda^3)+O(\epsilon).
  \end{equation}
  Similar calculations give the following values for the $\tilde N'_{0j}$'s, up to summands of order
  $\epsilon$,
\begin{equation}\label{Ntildes}
\tilde N'_{02}=\tilde
N'_{05}=-3(\tau-\tilde\kappa)(\tau^2-\lambda^2)+2(\tau^3-\lambda^3),\;\;
\; \tilde N'_{03}=3\tilde\kappa\tau^2 ,\;\;\; \tilde
N'_{04}=3\lambda^2(\tilde\kappa-\mu).
 \end{equation}
It follows from (\ref{Ntildes}), (\ref{N01t}) and (\ref{Ipsitil})
\begin{equation}\label{Itil}
I_{\tilde\psi}=6\tilde\kappa(2\tau^2-\lambda^2)-3(\tau^3-2\tau\lambda^2+\lambda^3).
\end{equation}
After (\ref{Kaptilde}) we obtain
$$I_{\tilde\psi}=-3I_{\psi}.$$

In the definition of $M$ the variables  $z_1,z_2,z_5$ play the
same role. However we can consider the following $S^1$ action on
$M$
\begin{equation}\label{DefHatps}
\Hat \psi_t[z]=[z_1,z_2,z_3,e^{2\pi it}z_4,z_5].
\end{equation}
Its Hamiltonian is $\Hat f=\pi\rho^2_4-\Hat \kappa,$ with
\begin{equation}\label{Kappah}
\Hat\kappa=\frac{1}{4}
\Big(\frac{\lambda^4-4\lambda\tau^3+3\tau^4}{\tau^3-\lambda^3}\Big).
\end{equation}

The corresponding $\Hat N'_{0j}$ have the following values up
summand of order $\epsilon$
\begin{equation}\label{Iha}
\Hat N'_{01}=\Hat N'_{02}=\Hat
N'_{05}=3(\lambda+\Hat\kappa)(\tau^2-\lambda^2)-2(\tau^3-\lambda^3),\,
\Hat N'_{03}=3\tau^2(\Hat\kappa-\mu),\, \Hat
N'_{04}=3\Hat\kappa\lambda^2.
\end{equation}
From the preceding formulae one deduces
$$I_{\Hat\psi}=-I_{\tilde\psi}=3I_{\psi}.$$

\begin{Thm}\label{FThm} Let $M$ be the toric manifold defined by (\ref{DefM})
and (\ref{ActD}). If $\psi$, $\tilde\psi$ and $\Hat\psi$ are the
Hamiltonian loops in $M$ defined by (\ref{DEfpsi}),
(\ref{Defetildepsi}) and (\ref{DefHatps}) respectively, then
$$I_{\Hat\psi}=-I_{\tilde\psi}=3I_{\psi},$$
with
$$I_{\psi}=
\frac{\lambda^2(-3\tau^4+8\tau^3\lambda-6\tau^2\lambda^2
+\lambda^4)}{2(\tau^3-\lambda^3)},$$
  $\lambda$ being $\lambda:=\tau-\mu$.
\end{Thm}

\begin{Cor} Let $(M,\omega)$ be the toric manifold one point blow up of ${\Bbb
C}P^2$, then $\pi_1(\text{Ham}(M,\omega))$ contains an infinite
cyclic subgroup.
\end{Cor}
\begin{proof}
By Proposition \ref{Prononv} $I_{\psi}\ne 0$. As $I$ is a group
homomorphism then the class
$[\psi^l]\in\pi_1(\text{Ham}(M,\omega))$ does not vanish, for all
$l\in{\Bbb Z}\setminus \{0\}$.
\end{proof}

\section{Hamiltonian group of toric manifolds}\label{blowup}
In this Section we generalize the calculations carried out in
Section 2 for the $6$-manifold one point blow up of ${\Bbb C}P^3$
to a general toric manifold.

 Let ${\Bbb T}$ be the torus
$(S^1)^r$, and ${\frak t}={\Bbb R}\oplus\dots\oplus{\Bbb R}$ its
Lie algebra. Given $w_{j}\in{\Bbb Z}^r$, with $j=1,\dots,m$
 and $\tau\in{\Bbb R}^r$ we put
 \begin{equation}\label{DefiM}
 M=\{z\in {\Bbb C}^m\,:\,\pi\sum_{j=1}
^m|z_{j}|^2w_{j}=\tau\}/{\Bbb T},
 \end{equation}
 where the relation defined by ${\Bbb T}$ is
 \begin{equation}\label{DefiAc}(z_j)\simeq(z'_j)\;\;\text{iff
there is}\;\,\xi\in {\frak t}\;\, \text{such that}\;\,
z'_j=z_je^{2\pi i\langle w_j,\xi \rangle}\; \,\text{for}\;\,
j=1,\dots, m.
 \end{equation}

We will assume that there is an open half space in ${\Bbb R}^r$
which contains all the vectors $w_j$ and that $\{w_j\}_j$ span
${\Bbb R}^r$. We also assume that $\tau$ ia a regular value of the
map
$$z\in{\Bbb C}^m\mapsto \pi\sum_{j=1}^m|z_j|^2w_j\in{\Bbb R}^r.$$
Then $M$ is a closed toric manifold of dimension $n:=2(m-r)$
\cite{Mc-S1}.

When $0\ne z_a\in{\Bbb C}$, we write $z_a=\rho_ae^{i\theta_a}$.
The standard symplectic form on  ${\Bbb C}^m$ gives rise to the
symplectic structure $\omega$ on $M$. On
 $$\{[z]\in M\,:\,z_j\ne 0\;\text{for all}\; j \}$$
 $\omega$ can be written as
 $$\omega=\sum_{i=1}^n d\Big(\frac{\rho_{ai}^2}{2}
 \Big)\wedge d\,\varphi_{ai},$$
 with $\varphi_{ai}$ a linear combination of the $\theta_c$'s.

 Given $0<\epsilon<<1$, we set
 $$B_0=\{[z]\in M\,:\,|z_j|>\epsilon \;\,\text{for all}\; j\}$$
$$B_k=\{[z]\in M\,:\,|z_k|<\epsilon,\;|z_j|>\epsilon \;\,\text{for
all}\; j\ne k\},$$
 as in Section 2.
 On $B_0$ we will consider the Darboux coordinates
 $$\{\frac{\rho_{ai}^2}{2},\,\varphi_{ai} \}_{i=1,\dots,n}.$$
 Given $k\in\{1,\dots,m\}$ we write $\omega$ in the form
 $$\omega=d\Big(\frac{\rho_k^2}{2}\Big)
\wedge d\varphi_k+\sum_{i=1}^{n-1}d\Big(\frac{\rho_{ki}^2}{2}\Big)
\wedge d\varphi_{ki},$$
 where $\varphi_k$ and $\varphi_{ki}$ are
linear combinations of the $\theta_c$'s. Then we consider on $B_k$
the following Darboux coordinates
 $$\{x_k,y_k, \frac{\rho_{ki}^2}{2}, \varphi_{ki} \}_{i=1,\dots,n-1},$$
 with $x_k+iy_k:=\rho_ke^{i\varphi_k}.$

 We denote by $\psi_t$ the map
 $$\psi_t:[z]\in M\mapsto[e^{2\pi it}z_1,z_2,\dots,z_m]\in M.$$
 $\{\psi_t\,:\, t\in[0,1]\}$ is a loop in $\text{Ham}(M)$. By repeating  the arguments of Section 2 one obtains
 $$I_{\psi}=\sum_{k=1}^mN'_{0k}+O(\epsilon),$$
 where
 $$N'_{0k}=\frac{ni}{2\pi}\int_{A'_{0k}}fd\,\log
 r_{0k}\wedge\omega^{n-1},$$
$$A'_{0k}=\{[z]\in M\,:\,|z_k|=\epsilon,\,|z_j|>\epsilon\; \text{for
all}\; j\ne k \},$$
 and $f=\pi\rho_1^2-\kappa$, with
 $$\int_M\pi\rho_1^2\omega^n=\kappa\int_M\omega^n.$$

 As in Section  2, on $A'_{0k}$ the exterior derivative  $d\,\log r_{0k}=ih'(\varphi_k)d\varphi_k,$
 where $h=h(\varphi_k)$ is a function such that
 $h(0)=2\pi,\,h(2\pi)=0$. Then
 $$N'_{0k}=-n\int_{\{ [z]\,:\,z_k=0
 \}}f\omega^{n-1}+O(\epsilon).$$
 Since $I_{\psi}$ is independent of $\epsilon$, we obtain
 \begin{equation}\label{Ipsito}
  I_{\psi}=-n \sum_{k=1}^m\Big(\int_{\{ [z]\,:\,z_k=0
 \}}(\pi\rho_1^2-\kappa)\omega^{n-1} \Big).
 \end{equation}
This formula together with the fact that $I$ is a group
homomorphism give the following Theorem
\begin{Thm}\label{Thmtoric} Let $(M,\omega)$ be the toric manifold defined
by (\ref{DefiM}) and (\ref{DefiAc}). If
$$\sum_{k=1}^m\Big(\int_{\{ [z]\,:\,z_k=0
 \}}(\pi\rho_1^2-\kappa)\omega^{n-1} \Big)\ne 0,$$
 then $\pi_1(\text{Ham}(M,\omega))$ contains an infinite cyclic
 subgroup.
\end{Thm}

\smallskip

{\it Examples.} We will  check the above result  calculating
$I_{\psi}$ by (\ref{Ipsito}) in two particular cases: When the
manifold is ${\Bbb C}P^1$ and when is ${\Bbb C}P^2$.

\smallskip

For
$${\Bbb C}P^1=\{[z_1,z_2]\,:\, |z_1|^2+|z_2|^2=\tau/\pi  \}/S^1$$
and $\psi_t[z_1,z_2]=[e^{2\pi it}z_1,z_2]$, the normalized
Hamiltonian is $f=\pi\rho_1^2-\tau/2$, that is, $\kappa=\tau/2.$
In this case (\ref{Ipsito}) reduces to $I_{\psi}=-\tau+2\kappa=0.$
This is compatible with the fact that $\pi_1(\text{Ham}({\Bbb
C}P^1))={\Bbb Z}/2{\Bbb Z}.$

\smallskip

For
$${\Bbb C}P^2=\{[z_1,z_2,z_3]\,:\, |z_1|^2+|z_2|^2+|z_3|^2=\tau/\pi  \}/S^1,$$
the Hamiltonian is $f=\pi\rho_1^2-\tau/3.$
 Moreover for for $k\in\{1,2,3\}$
 $$\int_{\{[z]\,:\,z_k=0  \}}\omega=\tau.$$
On the other hand, for $k=2,3$
$$\int_{\{[z]\,:\,z_k=0  \}}\pi\rho_1^2\omega=\tau^2/2.$$
After (\ref{Ipsito})  $I_{\psi}=-2(\tau^2-3\kappa\tau)=0$. This
result is consistent with the finiteness of
 $\pi_1(\text{Ham}({\Bbb
C}P^2))$, since   $\text{Ham}({\Bbb C}P^2)$ has the homotopy type
of $PU(3)$ \cite{Gr}.


\end{document}